\documentclass[12pt]{article}
\usepackage{maple2e,amsmath,amsopn,amssymb,amsthm,amsfonts}
\usepackage[T2A]{fontenc}
\usepackage[cp1251]{inputenc}
\usepackage{longtable}

\newcommand\g{{\mathfrak g}}

\newcommand{\R}{\mathbb{R}}
\newcommand{\Co}{\mathbb{C}}

\begin{document}

{\bf \large
\centerline{N.~K.~Smolentsev}

\vspace{3mm}
\centerline{Left-invariant almost para-complex structures}
\centerline{on six-dimensional nilpotent Lie groups}}
\vspace{3mm}

\begin{abstract}
There are five six-dimensional nilpotent Lie groups $G$, which do not admit neither symplectic, nor complex structures and, therefore, can be neither almost pseudo-K\"{a}hler, nor almost Hermitian. In this paper, we study precisely these Lie groups. The aim of the paper is to define new left-invariant geometric structures on the Lie groups under consideration that compensate, in some sense, the absence of symplectic and complex structures.
Weakening the closedness requirement of left-invariant 2-forms $\omega$ on the Lie groups, non-degenerated 2-forms $\omega$ are obtained, whose exterior differential $d\omega$ is also non-degenerated in Hitchin sense \cite{Hitchin}. Therefore, the Hitchin’s operator $K_{d\omega}$ is defined for the 3-form $d\omega$. It is shown that $K_{d\omega}$ defines an almost complex or almost para-complex structure for $G$ and the pair $(\omega, d\omega)$ defines pseudo-Riemannian metrics of signature (2.4) or (3.3), for which the Ricci operator is diagonal and has two eigenvalues differing in sign.
\end{abstract}

\section{Introduction} \label{Introduction}
Left-invariant K\"{a}hler structure on Lie group $G$ is a triple $(g, \omega, J)$ consisting of a  left-invariant Riemannian metric $g$, left-invariant symplectic form $\omega$ and orthogonal left-invariant complex structure $J$, where $g(X,Y) = \omega(X,JY)$ for any left-invariant vectors fields $X$ and $Y$ on $G$.
Therefore, such a structure on $G$ can be given by a pair $(\omega, J)$, where $\omega$ is a symplectic form, and $J$ is a complex structure compatible with $\omega$, that is, such that $\omega(JX,JY) = \omega(X,Y)$. If $\omega(X,JX) > 0$, $\forall X > 0$, it is K\"{a}hler metrics.
If the positivity condition is not satisfied, then $g(X,Y) = \omega(X,JY)$ is a pseudo-Riemannian metric and then $(g, \omega, J)$ is called a pseudo-K\"{a}hler structure on the Lie group $G$.
Classification of real six-dimensional nilpotent Lie algebras admitting invariant complex structures were obtained in \cite{Sal-1}.
Classification of symplectic structures on six-dimensional nilpotent Lie algebras was obtained in \cite{Goze-Khakim-Med}.
Out of 34 classes of isomorphic simply connected six-dimensional nilpotent Lie groups, only 26 admit left-invariant symplectic structures.
Condition of existence of left-invariant positively definite metric on Lie group $G$ applies restrictions to the structure of its Lie algebra $\g$. For example, it was shown in \cite{BG} that such a Lie algebra can not be nilpotent except for the abelian case.
Although nilpotent Lie groups and nilmanifolds (except for torus) do not admit K\"{a}hler left-invariant metrics, on such manifolds left-invariant pseudo-Riemannian K\"{a}hler metrics may exist.
It was shown in \cite{CFU2} that 14 classes of symplectic six-dimensional nilpotent Lie groups admit compatible complex structures and, therefore, define pseudo-K\"{a}hler metrics. A more complete study of the properties of the curvature of such pseudo-K\"{a}hler and almost pseudo-K\"{a}hler structures was carried out in \cite{Smolen-11, Smolen-13}. In \cite{Conti-18}, it was shown that on a nilpotent Lie algebra of dimension up to six, Einstein metrics are Ricci-flat.

As mentioned before, 26 out of 34 classes of six-dimensional nilpotent Lie groups admit left-invariant symplectic structures. Out of last 8 classes of non-symplectic Lie groups, 5 Lie groups $G_i$ do not also admit complex structures \cite{Sal-1}, their Lie algebras $\g_i$ are shown below:
$$
\begin{array}{l}
  \g_{1} = (0, 0, 12, 13, 14+23, 34 -25), \\
  \g_{2} = (0, 0, 12, 13, 14, 34 -25), \\
  \g_{3} = (0, 0, 0, 12, 13, 14+35),\\
  \g_{4} = (0, 0, 0, 12, 23, 14+35), \\
  \g_{5} = (0, 0, 0, 0, 12, 15+34).
\end{array}
$$

In this paper we study precisely these Lie groups. The aim of the paper is to define new left-invariant geometric structures on the Lie groups under consideration that compensate, in some sense, the failure of symplectic and complex structures.
It is shown that on all such Lie groups $G_i$ any left-invariant closed 2-form $\omega$ is degenerated. There are natural ways to weaken the closedness requirement of $\omega$ to preserve non-degeneracy $\omega$, in ways that 3-form $d\omega$ is also non-degenerated and property $\omega\wedge d\omega = 0$ is satisfied.
Hitchin's operator $K_{d\omega}$ corresponding to non-degenerated 3-form $d\omega$, can define either almost complex structure, or almost para-complex, depending on the chosen $\omega$.
Associated metric $g(X,Y) = \omega(X,J_{d\omega}Y)$ is pseudo-Riemannian of signatures (3,3) or (2,4). The structural group is reduced to $SL(3,\R)$ in case of signature (3,3) and to $SU(1,2)$ in case of signature (2,4). On the groups $G_2 - G_5$, these metrics have non-zero scalar curvature and the Ricci operator with two distinct sign eigenvalues. Pseudo-Riemannian metrics on the group $G_1$ are not Ricci-flat, but the scalar curvature may be zero for some values of the parameters.
An explicit form of these metrics is presented.
As a result, we obtain a compatible couple $(\omega,\Omega)$, where $\Omega = d\omega$. We present an explicit form of pseudo almost Hermitian half-flat and para-complex half-flat structures.

For any nilpotent Lie group $G$ with rational structure constants there exists a discrete subgroup $\Gamma$ such that $M = \Gamma\backslash G$ is a compact manifold called a nilmanifold. Therefore, all the results hold for the corresponding six-dimensional compact nilmanifolds.

All calculations were made in the Maple system according to the usual formulas for the geometry of left-invariant structures.

\section{Preliminaries} \label{Preface}
Let $G$ be a real Lie group of dimension $m$ and $\g$ be its Lie algebra. Lower central series of Lie algebra $\g$ is decreasing sequence of ideals $C^0\g$, $C^1\g$, \dots, being defined inductively: $C^0\g = \g$, $C^{k+1}\g = [\g, C^k\g]$.
Lie algebra $\g$ is called nilpotent, if $C^k\g = 0$ for some $k$.
In this case, the minimum length of lower central series is called class (or step) of nilpotency.
In other words, the Lie algebra class is equal to $s$, if $C^s\g = 0$ and $C^{s-1}\g \neq 0$.
In this case, $C^{s-1}\g$ lies in the center $\mathcal{Z}(\g)$ of the Lie algebra $\g$.
The increasing central sequence $\{g_l\}$ was defined for nilpotent $s$-step Lie algebra,
$$
\g_0 = {0} \subset \g_1 \subset \g_2 \subset \dots \subset \g_{s-1} \subset \g_s = \g,
$$
where the  ideals $g_l$ were defined inductively by the rule:
$$
g_l = \{X \in \g | [X, \g] \subseteq \g_{l-1}\}, \ l \ge 1.
$$
Particularly, $\g_1$ is the center of Lie algebra. One can see from this sequence that nilpotency property is equivalent to existence of basis $\{e_1, \dots, e_m\}$ of the Lie algebra $\g$, for which
$$
[e_i,e_j] =\sum_{k>i,j}C^k_{ij}e_k,\ 1\le i< j \le m.
$$
Nilpotency is also equivalent to the existence of basis $\{e^1, \dots, e^m\}$ of left-invariant 1-forms on $G$ such that
$$
de^i \in \Lambda^2\{e^1, \dots, e^{i-1}\},\quad  1\le i \le m,
$$
where the right side is considered to be zero for $i = 1$.
As is known, the exterior differential of a left-invariant 1-form is expressed through the structural constants of a Lie algebra by the formula \cite{KN}:
$$
de^k =-\sum_{i<j}C^k_{ij}e^i\wedge e^j,
$$
where $\{e^1, \dots, e^m\}$ is the dual basis in $\g^*$ to $\{e_1, \dots, e_m\}$.
Therefore 
the structure of Lie algebra is given either by specifying nonzero Lie brackets or by differentials of basis left-invariant 1-forms. The Lie algebra $\g$ is often defined as an $m$-tuple based on a sequence of differentials $(0, 0, de^3, \dots, de^m)$ of basic 1-forms, in the notation $ij$ is used instead of $e^{ij} = -e^i\wedge e^j$.
For example, notation $(0,0,0,0,12,14+23)$ denotes Lie algebra with structural equations: $de^1 = de^2 = de^3 = 0$, $de^4 = 0$, $de^5 = -e^1\wedge e^2$ and $de^6 =-e^1\wedge e^4= -e^2\wedge e^3$.
Equivalently, the basis $\{e_1,\dots,e_6\}$ of $\g$ satisfies $[e_1,e_2]=e_5,\,[e_1,e_4]=[e_2,e_3]=e_6$.

Left-invariant symplectic structure on Lie group $G$ is a left-invariant closed 2-form $\omega$  of the maximal rank. It is given by 2-form $\omega$ of the maximal rank on Lie algebra $\g$. Closedness of the 2-form is equivalent to condition
$$	
\omega([X,Y],Z) -\omega([X,Z],Y) + \omega([Y,Z],X) = 0,\ \ \forall X,Y,Z \in \g.
$$
In this case, Lie algebra $\g$ and group $G$ will be called symplectic ones.

Left-invariant almost complex structure on Lie group $G$ is left-invariant field of endomorphisms $J: TG \to TG$ of tangent bundle $TG$, having the property $J^2 = -Id$.
Since $J$ is defined by linear operator $J$ on Lie algebra $\g = T_eG$, we will say that $J$ is a left-invariant almost complex structure on Lie algebra $\g$.
In order for the almost complex structure $J$ to define a complex structure on the Lie group $G$, it is necessary and sufficient (according to the Newlender-Nirenberg theorem \cite{KN}) that the Nijenhuis tensor vanishes:
$$
[JX,JY] -[X,Y] -J[JX,Y] -J[X,JY] = 0,\  \forall X,Y \in g.
$$
For the left-invariant complex structure on Lie group $G$ left shifts $L_a:G\to G$, $a \in G$  are holomorphic.

Left-invariant K\"{a}hler structure on Lie group $G$ is a triple $(g, \omega, J)$ consisting of a  left-invariant Riemannian metric $g$, left-invariant symplectic form $\omega$ and orthogonal left-invariant complex structure $J$, where $g(X,Y) = \omega(X,JY)$, $\forall X,Y \in \g$.
Therefore, such a structure on Lie group $G$ can be specified by a couple $(\omega,J)$, where $\omega$  is a symplectic form, and $J$ is a complex structure being compatible with $\omega$, i.e. such that $\omega(JX,JY) = \omega(X,Y)$, $\forall X,Y \in \g$. If $\omega(X,JY) > 0$, $\forall X \neq 0$, then it is K\"{a}hler metric $g(X,Y) = \omega(X,JY)$. But if the positivity condition is not fulfilled, then $g(X,Y)$ is pseudo-Riemannian metric and then $(g,\omega, J)$ is called pseudo-K\"{a}hler structure on Lie group $G$. In further, (pseudo)K\"{a}hler structure will be specified by pair $(\omega,J)$ of compatible left-invariant complex and symplectic structures.
It follows from left-invariance that (pseudo)K\"{a}hler structure $(g,\omega,J)$ can be given by the values $J, \omega$ and $g$  on Lie algebra $\g$ of the Lie group $G$. In this case $(g,\omega,J,\g)$ is called \emph{pseudo-K\"{a}hler Lie algebra}.

\emph{Almost para-complex structure} on $2n$-dimensional manifold $M$ is a field $P$ of endomorphisms of the tangent bundle $TM$ such that $P^2 = Id$, where ranks of eigen-distributions $T^\pm M:= \rm{ker}(Id\pm P)$ are equal.
Almost para-complex structure $P$ is called integrable if distributions $T^\pm M$ are involutive. In this case, $P$ is called \emph{para-complex structure}.
A manifold $M$ supplied by (almost) para-complex structure $P$, is called (almost) para-complex manifold. The Nijenhuis tensor $N_P$ of almost para-complex structure $P$ is defined by equation
$$
N_P(X,Y) = [X,Y] + [PX,PY] -P[PX,Y] -P[X,PY],
$$
for all vector fields $X,Y$ on $M$.
As in the case with complex structure, para-complex one $P$ is integrable if and only if $N_P = 0$.

Para-K\"{a}hler manifold can be defined as pseudo-Riemannian manifold $(M,g)$ with skew-symmetric para-complex structure $P$, that is parallel with respect to the Levi-Civita connection.
If $(g,P)$ is a para-K\"{a}hler structure on $M$, then $\omega = g\cdot P$ is symplectic structure, and eigen-distributions $T^\pm M$, corresponding to eigen-values $\pm 1$ of field $P$, represent two integrable $\omega$-Lagrangian distributions. Therefore the para-K\"{a}hler structure can be identified with bi-Lagrangian structure $(\omega,T^\pm M)$, where $\omega$ is a symplectic structure, and $T^\pm M$  are the integrable Lagrangian distributions. 	
In \cite{Aleks} presents a review of the theory para-complex structures, and the invariant para-complex and para-K\"{a}hler structures on homogeneous spaces of semi-simple Lie groups are considered in detail. It is shown that every invariant para-K\"{a}hler structure $P$ on $M = G/H$ defines a unique para-K\"{a}hler Einstein structure $(g,P)$ with given non-zero scalar curvature.

Since the 2-form $\omega$ is not closed, it is possible to consider the 3-form $d\omega$. In \cite{Hitchin} Hitchin had defined the concept of non-degeneracy (stability) for 3-forms $\Omega$ and built a linear operator $K_{\Omega}$, whose square is proportional to identity operator $Id$. Recall the basic Hitchin's constructions.

Let $V$ be a six-dimensional real vector space, $\mu$ be a volume form on $V$, and $\Lambda^3 V^*$ be the 20-dimensional linear space of skew-symmetric 3-forms on $V$.
We shall take interior product $\iota_X \Omega \in \Lambda^2 V^*$ for the form $\Omega \in \Lambda^3 V^*$ and vector $X \in V$. Then $\iota_X \Omega\wedge\Omega \in \Lambda^5 V^*$. Natural pairing by the exterior product $V^* \otimes \Lambda^5 V^* \to \Lambda^6 V^* \cong \R \mu$ defines the isomorphism $A: \Lambda^5 V^* \to V$. Using $\Lambda^5 V^* \cong V$ we define linear map $K_\Omega: V \to V$ as
$$
K_\Omega(X) = A\left(\iota_X \Omega\wedge\Omega\right).
$$

In other words, $\iota_{K_\Omega(X)} \mu = \iota_X\Omega\wedge\Omega$. Define $\lambda(\Omega) \in \R$ as a trace of the square of $K_\Omega$:
$$
\lambda(\Omega) = \frac 16{\rm tr}\,K_\Omega^2.
$$

The form $\Omega$ is called non-degenerated (or stable) if $\lambda(\Omega)\neq 0$.

It is shown in \cite{Hitchin} that if $\lambda(\Omega)\neq 0$, then
\begin{itemize}
  \item $\lambda(\Omega)> 0$ if and only if $\Omega = \alpha + \beta$, where $\alpha$, $\beta$ are real decomposable 3-forms and $\alpha\wedge \beta\neq 0$;
  \item $\lambda(\Omega)< 0$ if and only if $\Omega = \alpha + \overline{\alpha}$, where $\alpha \in \Lambda^3(V^*\otimes \Co)$ is complex decomposable 3-form and $\alpha\wedge\overline{\alpha} \neq 0$.
\end{itemize}

It follows that if $\Omega$ is real and $\lambda(\Omega)> 0$, then it lies in $GL(V)$-orbit of form $\varphi = \theta^1\wedge\theta^2\wedge\theta^3 + \theta^4\wedge\theta^5\wedge\theta^6$ for some basis $\theta^1,\dots \theta^6$ in $V^*$, and if $\lambda(\Omega)<0$, then it lies in orbit of form $\varphi = \alpha + \overline{\alpha}$, where $\alpha = (\theta^1 +i\theta^2)\wedge (\theta^3 +i\theta^4)\wedge(\theta^5 +i\theta^6)$.

Then the real 20-dimensional vector space $\Lambda^3(V^*)$ contains invariant quadratic hypersurface $\lambda(\Omega)= 0$ dividing the $\Lambda^3(V^*)$ to 2 open sets: $\lambda(\Omega)>0$ and $\lambda(\Omega)< 0$.
The component of the unit of the stabilizer of the 3-form lying in the first set is conjugate to the group $SL(3,\R)\times SL(3,\R)$, and in the other case to the group $SL(3,\Co)$.
The linear transformation of $K_\Omega$ has \cite{Hitchin} the following properties: $\rm{tr}\,K_\Omega = 0$ and $K_\Omega^2= \lambda(\Omega)\, Id$.
In the case $\lambda(\Omega)<0$, the real 3-form $\Omega$ defines the structure $J_\Omega$ of complex vector space on real vector space $V$ as follows:
$$
J_\Omega \frac{1}{\sqrt{-\lambda(\Omega)}}K_\Omega.
$$
But if $\lambda(\Omega)>0$, 3-form $\Omega$ defines the para-complex structure $J_\Omega$, i.e. $J_\Omega^2 = 1$, $J_\Omega\neq 1$  on real vector space $V$ by similar formula:
$$
J_\Omega \frac{1}{\sqrt{\lambda(\Omega)}}K_\Omega.
$$
Recall that the structure of almost a product is called para-complex, if eigen-subspaces have the same dimension.

The elements of $GL(V)$-orbits of 3-form $\Omega$, corresponding to $\lambda(\Omega)>0$, have stabilizer  $SL(3,\R)\times SL(3,\R)$ in $GL^+(V)$ and $J_\Omega$ is para-complex structure, i.e. $J_\Omega^2 = 1$, $J_\Omega\neq 1$.
The elements of orbit corresponding to $\lambda(\Omega)<0$, have stabilizer $SL(3,\Co)$ in $GL^+(V)$ and $J_\Omega$ is almost complex structure, i.e. $J_\Omega^2 = -1$.
In both cases, dual to $\Omega$ form is defined by formula $\widehat{\Omega} = J_\Omega^*\Omega$.
If $\lambda(\Omega)>0$ and $\Omega = \alpha + \beta$, then  $\widehat{\Omega} = \alpha - \beta$.
But if $\lambda(\Omega)<0$ and $\Omega = \alpha + \overline{\alpha}$, then $\widehat{\Omega} = i(\overline{\alpha}-\alpha)$.
Note that $\Omega^{\wedge\wedge} = -\Omega$ in both cases and $J_{\widehat{\Omega}} =-\varepsilon J_\Omega$, where $\varepsilon$ is the sign of $\lambda(\Omega)$.
The additional 3-form $\widehat{\Omega}$ has a defining property: if $\lambda(\Omega)>0$, then $\Psi =\Omega+\widehat{\Omega}$ is decomposable, and if $\lambda(\Omega)<0$, then complex form $\Psi =\Omega +i\widehat{\Omega}$ is decomposable.

 The pair $(\omega,\Omega)\in \Lambda^2(V^*)\times \Lambda^3(V^*)$ of non-degenerated forms is called \emph{compatible} if $\omega\wedge \Omega =0$ (or, equivalently, $\widehat{\Omega}\wedge \omega =0$), and it is called normalized, if $\widehat{\Omega}\wedge \Omega = 2\omega^3/3$.

 Each compatible pair $(\omega,\Omega)$ uniquely defines $\varepsilon$-complex structure $J_\omega$ (i.e. $J_\omega^2 =\varepsilon$), scalar product $g_{(\omega,\Omega)}(X,Y) = \omega(X,J_\omega Y)$ (signatures (3,3) for $\varepsilon =1$ and signatures (2,4) or (4,2) for $\varepsilon =-1$), and (para-)complex volume form $\Psi =\Omega +i_{\varepsilon}\widehat{\Omega}$ of type (3,0) with respect to $J_\omega$ (where $i_{\varepsilon}$ is a complex or para-complex imaginary unit).
 In addition, the stabilizer of $(\omega,\Omega)$ pair is $SU(p,q)$ for $\varepsilon =-1$ and $SL(3,\R)\subset SO(3,3)$ for $\varepsilon =1$. Therefore, $(\omega,\Omega)$ pair for $\varepsilon =-1$ defines pseudo almost Hermitian structure. But if $\varepsilon =1$, it defines almost para-Hermitian structure. Such structures are also called special almost $\varepsilon$-Hermitian.

 Also recall that $SU(3)$ structure on real six-dimensional almost Hermitian manifold $(M, g, J, \omega)$ is specified by $(3,0)$-form $\Psi$.
 Almost Hermitian 6-manifold is called half-flat \cite{Conti} if it admits a reduction to $SU(3)$, for which $d \Re(\Psi) = 0$ and $\omega\wedge d\omega =0$.

 In the case of pseudo-Riemannian manifold, each compatible pair $(\omega,\Omega)$ uniquely defines the reduction to $SU(1,2)$ for $\varepsilon =-1$ and to $SL(3,\R)\subset SO(3,3)$ for $\varepsilon =1$. Therefore, 6-manifold with the $(\omega,\Omega)$ pair possessing the properties $d\Omega = 0$ and $\omega\wedge d\omega =0$, will be called \emph{half-flat pseudo almost Hermitian} if it admits the reduction to $SU(1,2)$ or \emph{half-flat almost para-complex} if it admits the reduction to $SL(3,\R)\subset SO(3,3)$.

 \textbf{Remark.} In this work, we assume that exterior product and exterior differential are defined without normalizing constant. In particular, then $dx\wedge dy = dx\otimes dy -dy\otimes dx$ and $d\omega(X,Y) = X\omega(Y) -Y\omega(X) -\omega([X,Y])$. Let $\nabla$ be the Levi-Civita connectivity corresponding to left-invariant (pseudo)Riemannian metric $g$. It is defined by six-membered formula \cite{KN}, which becomes the following form for left-invariant vector fields $X,Y,Z$ on Lie group: $2g(\nabla_XY,Z) = g([X,Y],Z) + g([Z,X],Y) + g(X,[Z,Y])$. If $R(X,Y) = [\nabla_X, \nabla_Y] -\nabla_{[X,Y]}$ is a curvature tensor, then Ricci tensor $Ric(X,Y)$ for (pseudo)Riemannian metric $g$ is defined as a construction of a curvature tensor over the first and fourth (upper) indices.

\section{Singular Lie groups} \label{Singular_Lie_G}
In this section Lie groups that do not admit neither symplectic, nor left-invariant complex structures will be considered. Such Lie groups will be called \emph{singular}. It will be shown that they admit non-degenerated left-invariant 2-forms, whose exterior differentials are non-degenerated.
In addition, they admit almost para-complex structures and pseudo-Riemannian metrics of signature (3.3), which have a diagonal Ricci operator with two eigenvalues, differing by the sign.

\subsection{Lie group $G_1$} \label{G_1}
Singular group $G_1$ that does not admit neither symplectic, nor complex structures. Non-zero commutation relations:
$$
[e_1,e_2] = e_3, [e_1,e_3] = e_4,, [e_1,e_4] = e_5, [e_2,e_3] = e_5, [e_3,e_4] = e_6, [e_2,e_5] = -e_6.
$$
Lie algebra $\g$ has ideals: $C^1\g = D^1\g = \R\{e_3,e_4,e_5,e_6\}$, $C^2\g = \R\{e_4,e_5,e_6\}$, $C^3\g = \R\{e_5,e_6\}$, $C^4\g = \mathcal{Z} = \R{e_6}$ is the center of Lie algebra. Filiform Lie algebra. Does not admit half-flat structures \cite{Conti}.

Let $\omega = a_{ij}e_i\wedge e_j$ be any left-invariant 2-form. For the general form $\omega$ Hitchin’s operator $K_\omega$ has quite complicated form and the following function $\lambda(d\omega)$:
$$
\lambda = 4\left(a_{16}a_{56}^2 + 4a_{35}a_{56}^2 + 4a_{36}^2a_{56} -4a_{36}a_{46})2 -4a_{45}a_{46}a_{56}\right)a_{56}+a_{46}^4.
$$
Thus, in general, 3-form $d\omega$ is non-degenerated.
It is easy to see that 2-form $\omega$ is closed if and only if $a_{16} = a_{26} = a_{36} = a_{35} = a_{45} = a_{46} = a_{56} = 0$ and $a_{34} = -a_{25}$, $a_{24} = a_{15}$.
However, such 2-form $\omega$ is degenerated.
There are several natural ways to weaken the closedness requirement of the 2-form $\omega$, so as not to lose the nondegeneracy of $\omega$ and $d\omega$.

\subsubsection{Option 1.}
In that case we will not suppose that coefficients $a_{46}$ and $a_{56}$, which essentially occur in the expression for function $\lambda(d\omega)$, are non-zero.
Moreover, we will suppose that $a_{56}\neq 0$.
Then the property $\omega\wedge d\omega = 0 $is fulfilled under condition $a_{15} = 0$, $a_{25} = 0$ and $a_{12}a_{56} =a_{13}a_{46}$, $a_{23} = -a_{14}$. 2-form $\omega$ is non-degenerated under condition that $a_{14}a_{56}\neq 0$ and $\omega$ and $d\omega$ become
$$
\omega = e^1\wedge(a_{13}a_{46}/a_{56}\, e^2 + a_{13}\, e^3 +a_{14}\,e^4) -a_{14}\, e^2\wedge e^3 + a_{46}\,e^4\wedge e^6 +a_{56} \,e^5\wedge e^6,
$$
$$
d\omega = -a_{46} e^{136} + a_{46} e^{245} -a_{56} e^{146} -a_{56} e^{236} +a_{56} e^{345}.
$$
The function $\lambda(d\omega)$ of the Hitchin's operator \cite{Hitchin} for 3-form $d\omega$ becomes $\lambda = a_{46}^4$. The Hitchin's operator $K_{d\omega}$ has a matrix
$$
K_{d\omega} = \left[ \begin {array}{cccccc}
-{{\it a_{46}}}^{2}&-2\,{\it a_{46}}\,{\it a_{56}}&-2\,{{\it a_{56}}}^{2}&0&0&0\\
\noalign{\medskip}0&{{\it a_{46}}}^{2}&2\,{\it a_{46}}\,{\it a_{56}}&2\,{{\it a_{56}}}^{2}&0&0\\ \noalign{\medskip}0&0& -{{\it a_{46}}}^{2}&-2\,{\it a_{46}}\,{\it a_{56}}&0&0\\
\noalign{\medskip}0&0&0&{{\it a_{46}}}^{2}&0&0\\
\noalign{\medskip}0&0&0&0&{{\it a_{46}}}^{2}&2\,{{\it a_{56}}}^{2}\\
\noalign{\medskip}0&0&0&0&0&-{{\it a_{46}}}^{2}
\end {array} \right]
$$

Determine the operator $P=K_{d\omega}/a_{46}^2$.
It defines left-invariant almost para-complex structure, $P^2 = Id$, having the property $\omega(PX,PY) = -\omega(X,Y)$.
Eigen-subspaces $E^\pm$ related to the eigen-values $\pm 1$ of operator $P$ are generated by the following vectors:
$$
E^+ = \{a_{56}^3\,e_1 -a_{56}a_{46}^2\,e_3 +a_{46}^3\,e_4,  -a_{56}\,e_1 +a_{46}\,e_2, e_5\},
$$
$$
E^- = \{\,e_1,  -a_{56}\,e_2 +a_{46}\,e_3,  -a_{56}^2\,e_5 +a_{46}^2\,e_6\}.
$$
It is easy to see that they are not closed relative to Lie bracket, so $P$ defines non-integrable almost a para-complex structure.

Define pseudo-Riemannian metric $g(X,Y) = \omega(X,PY)$ of signature (3,3).
It is given by:
$$
g =  \left[ \begin {array}{cccccc}
0&{\frac {{\it a_{13}}\,{\it a_{46}}}{{\it a_{56}}}}&{\it a_{13}}&{\it a_{14}}&0&0\\
\noalign{\medskip}{\frac {{\it a_{13}}\,{\it a_{46}}}{{\it a_{56}}}}&2\,{\it a_{13}}&{\frac {2\,{\it a_{13}}\,{\it a_{56}}+{\it a_{14}}\,{\it a_{46}}}{{\it a_{46}}}}&2\,{\frac {{\it a_{14}}\,{\it a_{56}}}{{\it a_{46}}}}&0&0\\ \noalign{\medskip}{\it a_{13}}&{\frac {2\,{\it a_{13}}\,{\it a_{56}}+{\it a_{14}}\,{\it a_{46}}}{{\it a_{46}}}}&2\,{\frac {{\it a_{56}}\, \left( {\it a_{13}}\,{\it a_{56}}+{\it a_{14}}\,{\it a_{46}} \right) }{{{\it a_{46}}
}^{2}}}&2\,{\frac {{\it a_{14}}\,{{\it a_{56}}}^{2}}{{{\it a_{46}}}^{2}}}&0&0\\
\noalign{\medskip}{\it a_{14}}&2\,{\frac {{\it a_{14}}\,{\it a_{56}}}{{\it a_{46}}}}&2\,{\frac {{\it a_{14}}\,{{\it a_{56}}}^{2}}{{{\it a_{46}}}^{2}}}&0&0&-{\it a_{46}}\\
\noalign{\medskip}0&0&0&0&0&-{\it a_{56}}\\
\noalign{\medskip}0&0&0&-{\it a_{46}}&-{\it a_{56}}&-2\,{\frac {{{\it a_{56}}}^{3}}{{{\it a_{46}}}^{2
}}}\end {array} \right]
$$

Direct calculations of curvature tensor in Maple system show that this metric is not Einsteinian and that it has scalar curvature
$$
R=\frac{8a_{13}a_{56}^7 -8a_{14}a_{46}a_{56}^6 -a_{46}^8}{a_{14}^2a_{46}^6a_{56}}.
$$
For $a_{13}: = (a_{46}^8 + 8 a_{56}^6 a_{46} a_{14})/(8 a_{56}^7)$ we get an example of a pseudo-Riemannian metric with zero scalar curvature and non-zero Ricci tensor.

\subsubsection{Option 2.}
Take the 2-form $\omega$ in the view $\omega = \omega_0 + \omega_C$, where $\omega_0$ is a closed 2-form and $\omega_C$ is a non-degenerated 2-form on the ideality $C^2\g = \R\{e_4,e_5,e_6\}$. We require that the 2-form $\omega$ should have the property $\omega\wedge d\omega = 0$:
$$
a_{25} = 0, a_{15} = 0, a_{12} a_{56} = a_{13}a_{46}, a_{56} a_{23} + a_{14}a_{56} = 0.
$$
Then $\omega$ is non-degenerated under condition $a_{14}a_{56} \neq 0$. The $\omega$ and $d\omega$ take the view:
$$
\omega = e^1\wedge\left(a_{13}a_{46}/a_{56} e^2 + a_{13} e^3 + a_{14} e^4\right) -a_{14} e^2\wedge e^3 + a_{45} e^4\wedge e^5 + a_{46} e^4\wedge e^6 + a_{56} e^5\wedge e^6,
$$
$$
d\omega = a_{45} e^{234} -a_{45} e^{135} -a_{46} e^{136} +a_{46} e^{245} -a_{56} e^{146} -a_{56} e^{236} +a_{56} e^{345}.
$$
In this case, $\lambda(d\omega)$ function is expressed by the $\lambda = a_{46}^4 -4a_{46}a_{45}a_{56}^2$ and it can take both positive and negative values.

\subsubsection{Case 1}
The function $\lambda(d\omega)$ takes the value $-1$ when $a_{45} = (a_{46}^4 + 1)/(4a_{46}a_{56}^2)$. Then the operator $J = K_{d\omega}$ defines almost complex structure compatible with $\omega$ and has the form:
$$
J= \left[ \begin {array}{cccccc}
-{{\it a_{46}}}^{2}&-2\,{\it a_{46}}\,{\it a_{56}}&-2\,{{\it a_{56}}}^{2}&0&0&0\\
\noalign{\medskip}{\frac {1+{{\it a_{46}}}^{4}}{{2\,\it a_{46}}\,{\it a_{56}}}}&{{\it a_{46}}}^{2}&2\,{\it a_{46}}\,{\it a_{56}}&2\,{{\it a_{56}}}^{2}&0&0\\
\noalign{\medskip}0&0&-{{\it a_{46}}}^{2}&-2\,{\it a_{46}}\,{\it a_{56}}&0&0\\
\noalign{\medskip}0&0&{\frac {1+{{\it a_{46}}}^{4}}{{2\, \it a_{46}}\,{\it a_{56}}}}&{{\it a_{46}}}^{2}&0&0\\
\noalign{\medskip}0&0&-{\frac {1+{{\it a_{46}}}^{4}}{{{2\, \it a_{56}}}^{2}}}&-{\frac {1 +{{\it a_{46}}}^{4}}{{2\, \it a_{46}}\,{\it a_{56}}}}&{{\it a_{46}}}^{2}&2\,{{\it a_{56}}}^{2}\\
\noalign{\medskip}0&0&{\frac {( 1+{{\it a_{46}}}^{4}) ^{2}}{{{8\, \it a_{46}}}^{2}{{\it a_{56}}}^{4}
}}&0&-{\frac {1+{{\it a_{46}}}^{4}}{{{2\, \it a_{56}}}^{2}}}&-{{\it a_{46}}}^{2}
\end {array} \right]
$$

Specify the associated pseudo-Riemannian metric by formula $g(X,Y) = \omega(X,JY)$ of signature (2,4). Direct calculations of curvature tensor in Maple system show that this metric is not Einsteinian and that it has scalar curvature
$$
R={\frac {8\,{{\it a_{56}}}^{7}{\it a_{13}}-8\,{{\it a_{56}}}^{6}{\it a_{46}}\,{\it
a_{14}}-1}{{{\it a_{14}}}^{2}{\it a_{56}}}}
$$
For $a_{13}: = (1+ 8 a_{56}^6 a_{46} a_{14})/(8 a_{56}^7)$ we get an example of a pseudo-Riemannian metric with zero scalar curvature and non-zero Ricci tensor.

\subsubsection{Case 2}
The function $\lambda(d\omega)$ takes the value $+1$ when $a_{45} = (a_{46}^4 -1)/(4a_{46}a_{56}^2)$. Then the operator $P = K_{d\omega}$ defines almost para-complex structure compatible with $\omega$ and $P$ has the same matrix, as the above almost complex structure $J$ has, where it is necessary to substitute $(a_{46}^4 -1)$ instead of $(a_{46}^4 +1)$. The associated metric $g(X,Y) = \omega(X,PY)$ is pseudo-Riemannian of signature (3,3); it is not the Einsteinian one and has the same scalar curvature as in the first case.

\textbf{Conclusions.} Any left-invariant closed 2-form $\omega$ on Lie group $G_1$ is degenerated. There are several ways to weaken the closedness requirement of $\omega$ to preserve non-degeneracy $\omega$, in ways that 3-form $d\omega$ is non-degenerated and the property $\omega\wedge d\omega = 0$ is fulfilled. Hitchin's operator $K_{d\omega}$ corresponding to non-degenerated 3-form $d\omega$, can define either almost complex structure, or almost para-complex, depending on the chosen $\omega$. Associated metric $g(X,Y) = \omega(X,J_{d\omega}Y)$ is pseudo-Riemannian of signature (2,4) or (3,3). As a result, we have obtained a compatible pair $(\omega,\Omega)$, where $\Omega = d\omega$. Therefore, the properties $d\Omega = 0$ and $\omega\wedge d\omega = 0$ are fulfilled in an obvious way. The (3,0)-form has a view of $\Psi = d\omega +i_{\varepsilon}\widehat{d\omega}$, where $i_{\varepsilon}$ is a complex or para-complex unit. Thus, half-flat pseudo almost Hermitian and half-flat para-complex structures were naturally defined on Lie group $G_1$.

\subsection{Lie group $G_2$} \label{G_2}
Singular group $G_2$ that does not admit neither symplectic, nor complex structures. Commutation relations
$$
[e_1,e_2] = e_3,	[e_1,e_3] = e_4,	[e_1,e_4] = e_5,	[e_3,e_4] = e_6,	[e_2,e_5] = -e_6.
$$
Lie algebra $\g$ has ideals: $C^1\g = D^1\g = \R\{e_3,e_4,e_5,e_6\}$, $C^2\g = \R\{e_4,e_5,e_6\}$, $C^3\g = \R\{e_5,e_6\}$, $C^4\g = \mathcal{Z} = \R\{e_6\}$ is the center of Lie algebra. Filiform Lie algebra. Does not admit half-flat structures \cite{Conti}.

Let $\omega = a_{ij}e_i\wedge e_j$ be any left-invariant non-degenerated 2-form. For such a generic 2-form the square of the Hitchin's operator \cite{Hitchin} for 3-form $d\omega$ has a diagonality: $K_{d\omega} = (a_{46}^2 -2a_{36}a_{56})^2 Id$.
Therefore, 3-form $d\omega$ is non-degenerated if $a_{46}^2 -2a_{36}a_{56}\neq 0$. The 2-form $\omega$ is closed only in the case when it has the form:
$$
\omega = e^1\wedge (a_{12} e^2 +a_{13} e^3 +a_{14} e^4 +a_{15} e^5) +e^2\wedge (a_{23} e^3 -a_{34} e^5) +a_{34} e^3\wedge e^4.
$$
Such 2-form $\omega$ is non-degenerated. In order to preserve the non-degeneracy of the $\omega$ and $d\omega$ at the minimal weakening of closedness property of $\omega$, two variants are possible: $a_{46}\neq 0$, or $a_{36}\neq 0$ and $a_{56}\neq 0$.
However, if $a_{56}\neq 0$, then simple calculations show that the property $\omega\wedge d\omega =0$ is incompatible with the non-degeneracy $\omega$.

Therefore, consider a case when $a_{46}\neq 0$. Then $K_{d\omega} = a_{46}^4 Id$. In addition, $\omega\wedge d\omega =0$ under condition that $a_{13} = 0$ and $a_{34} = 0$.
Then the 2-form $\omega$ is non-degenerated under condition $a_{23}a_{15}a_{46}\neq 0$, and the $\omega$ and $d\omega$ take the view:
$$
\omega = e^1\wedge (a_{12} e^2 +a_{14} e^4+ a_{15} e^5) +a_{23} e^2\wedge e^3 + a_{46} e^4\wedge e^6,
$$
$$
d\omega = a_{46}(-e^{136} + e^{245}).
$$
The operator $K_{d\omega}$ for the 3-form $d\omega$ has the diagonal form: $K_{d\omega} = {\rm diag}\{-a_{46}^2, a_{46}^2, -a_{46}^2, a_{46}^2, a_{46}^2, -a_{46}^2\}$.
Define the operator $P=K_{d\omega}/a_{46}^2$.
It defines left-invariant almost para-complex structure, $P^2 = Id$, having the property $\omega(PX,PY) = -\omega(X,Y)$.
Eigen-subspaces $E^\pm$ related to the eigen-values $\pm 1$ of operator $P$ are generated by the following vectors:
$$
E^+ = \{e_2, e_4, e_5 \}, E^- = \{e_1, e_3, e_6\}.
$$
It is easy to see that they are not closed relative to Lie bracket, so $P$ defines non-integrable almost a para-complex structure.

Define the pseudo-Riemannian metric $g(X,Y) =\omega(X,PY)$. It has a signature (3,3) and it is given by:
$$
g= \left[ \begin {array}{cccccc} 0&{\it a_{12}}&0&{\it a_{14}}&{\it a_{15}}&0\\
\noalign{\medskip}{\it a_{12}}&0&-{\it a_{23}}&0&0&0\\
\noalign{\medskip}0&-{\it a_{23}}&0&0&0&0\\
\noalign{\medskip}{\it a_{14}}&0&0&0&0&-{\it a_{46}}\\
\noalign{\medskip}{\it a_{15}}&0&0&0&0&0\\
\noalign{\medskip}0&0&0&-{\it a_{46}}&0&0\end {array} \right]
$$
Direct calculations in the Maple system show that for $a_{14} = 0$, this metric has a diagonal Ricci operator with two eigenvalues differing by the sign:
$$
RIC = \frac{a_{46}}{2a_{15}a_{23}}diag\{-1,-1,-1,+1,-1,+1\}.
$$

\subsection{Lie group $G_3$} \label{G_3}
Singular group $G_3$ that does not admit neither symplectic, nor complex structures. Commutation relations
$$
[e_1,e_2] = e_4, [e_1,e_3] = e_5, [e_1,e_4] = e_6,	[e_3,e_5] = e_6.
$$
Lie algebra $\g$ has ideals: $C^1\g = D \g = \R\{e_4,e_5,e_6\}$, $C^2\g =\R\{e_6\} = \mathcal{Z}$ is the center of Lie algebra Does admit half-flat structure \cite{Conti}.

Let $\omega = a_{ij}e_i\wedge e_j$  be any 2-form. The Hitchin's operator $K_{d\omega}$ for generic 2-form $\omega$ has a quite complicated view. Moreover, $K_{d\omega}^2 = a_{46}^4 Id$. For $\lambda = a_{46}^4 \neq 0$ the 3-form $d\omega$ is non-degenerated. The operator $P = K_{d\omega}/a_{46}^2$ defines the left-invariant almost para-complex structure on $\g$. The $\omega\wedge d\omega =0$ property is fulfilled under the following conditions:
$$
a_{12} a_{46} -a_{14} a_{26} -a_{23} a_{56} + a_{24} a_{16} +a_{25} a_{36} -a_{35} a_{26} =0,
$$
$$
a_{25} a_{46} -a_{24} a_{56} -a_{26} a_{45} = 0, \qquad	a_{35} a_{46} -a_{36} a_{45} -a_{34} a_{56} = 0.
$$

It is easy to see that the 2-form $\omega$ is closed only if
$$
\omega = e^1\wedge (a_{12}e^2+a_{13}e^3+a_{14}e^4+a_{15}e^5)+e^2\wedge a_{23}e^3+a_{24}e^4+a_{25}e^5)+e^3\wedge (a_{25}e^4+a_{35}e^5).
$$
In order to preserve the non-degeneracy of the $\omega$ and $d\omega$ at the minimal weakening of closedness property of $\omega$, consider the case when $a_{46} \neq 0$. The $\omega\wedge d\omega =0$ property is fulfilled under the condition $a_{12} = 0$, $a_{25} = 0$, $a_{35} = 0$. Thus, 2-form $\omega$ is non-degenerated if $a_{15}a_{23}a_{46} \neq 0$ and then we obtain:
$$
\omega = e^1\wedge (a_{13}e^3+a_{14}e^4+a_{15}e^5) + e^2\wedge (a_{23}e^3+a_{24}e^4) + a_{46} e^4\wedge e^6,
$$
$$
d\omega = -a_{46}(e^{126} + e^{345}). 	
$$

The operator $K_{d\omega}$ for the 3-form $d\omega$ has the diagonal view, $K_{d\omega} = {\rm diag}\{-a_{46}^2, -a_{46}^2, a_{46}^2, a_{46}^2, a_{46}^2, -a_{46}^2\}$.
Define the operator $P=K_{d\omega}/a_{46}^2$.
It defines left-invariant almost para-complex structure, $P^2 = Id$, having the property $\omega(PX,PY) = -\omega(X,Y)$.
Eigen-subspaces $E^\pm$ related to the eigen-values $\pm 1$ of operator $P$ are generated by the following vectors:
$$
E^+ = \{e_3, e_4, e_5 \},\qquad	E^- =\{e_1, e_2, e_6\}.
$$
It is easy to see that they are not closed relative to Lie bracket, so $P$ sets non-integrable almost a para-complex structure.

Define the pseudo-Riemannian metric $g(X,Y) = \omega(X,PY)$ of signature (3,3).
It is given by:
$$
g = \left[ \begin {array}{cccccc} 0&0&{\it a_{13}}&{\it a_{14}}&{\it a_{15}}&0\\
\noalign{\medskip}0&0&{\it a_{23}}&{\it a_{24}}&0&0\\
\noalign{\medskip}{\it a_{13}}&{\it a_{23}}&0&0&0&0\\
\noalign{\medskip}{\it a_{14}}&{\it a_{24}}&0&0&0&-{\it a_{46}}\\
\noalign{\medskip}{\it a_{15}}&0&0&0&0&0\\
\noalign{\medskip}0&0&0&-{\it a_{46}}&0&0\end {array} \right]
$$

Direct calculations in the Maple system show that for $a_{14} = a_{24} =0$, this metric has a diagonal Ricci operator with two eigenvalues differing by the sign:
$$
RIC = \frac{a_{46}}{2a_{15}a_{23}}diag\{-1,-1,-1,+1,-1,+1\}.
$$

\subsection{Lie group $G_4$} \label{G_4}
Singular group that does not admit neither symplectic, nor complex structures. Commutation relations:
$$
[e_1,e_2] = e_4, [e_2,e_3] = e_5, [e_1,e_4] = e_6,	[e_3,e_5] = e_6.
$$
Lie algebra $\g$ has ideals: $C^1\g = D^1\g = \R\{e_4,e_5,e_6\}$, $C^2\g =\R\{e_6\} = \mathcal{Z}$ is the center of Lie algebra. Does admit half-flat structure \cite{Conti}.

Let $\omega = a_{ij}e_i\wedge e_j$ be any 2-form. The Hitchin's operator $K_{d\omega}$ for generic 2-form $\omega$ has a quite complicated form. Moreover, $K_{d\omega}^2 = (a_{46}^2 -a_{56}^2)^2 Id$.
The $\omega\wedge d\omega =0$ property is fulfilled under the following conditions:
$$
-a_{12} a_{46} +a_{14} a_{26} -a_{16} a_{24} +a_{23} a_{56} -a_{25} a_{36} + a_{26} a_{35} = 0,
$$
$$
-a_{14}a_{56} +a_{15}a_{46} -a_{16}a_{45} = 0,\quad  a_{34} a_{56} -a_{35} a_{46} +a_{36} a_{45} = 0.
$$
It is easy to see that the 2-form $\omega$ is closed only if
$$
\omega = e^1\wedge (a_{12}e^2 +a_{13}e^3 +a_{14}e^4 +a_{15}e^5) +e^2\wedge (a_{23}e^3 +a_{24}e^4 +a_{25}e^5) +e^3\wedge (-a_{15}e^4 +a_{35}e^5).
$$
In order to preserve the non-degeneracy of the $\omega$ and $d\omega$ at the minimal weakening of closedness property of $\omega$, consider the case when $a_{46}\neq 0$ and $a_{56}\neq 0$.
The $\omega\wedge d\omega = 0$ property is fulfilled under the following conditions:
$$
-a_{12}a_{46} +a_{23}a_{56} = 0,\ -a_{14}a_{56} +a_{15}a_{46} = 0, \ -a_{15}a_{56} -a_{35}a_{46} =0
$$
and $d\omega$ is given by the sum of two decomposable 3-forms:
$$
d\omega = (a_{56} e^1 -a_{46} e^3)\wedge e^{45} +(-a_{46}e^1 +a_{56} e^3)\wedge e^{26}.
$$  	
The function $\lambda(d\omega)$ of the Hitchin's operator for 3-form $d\omega$ has the same view $\lambda = (a_{46}^2 -a_{56}^2)^2$.
And operator $K_{d\omega}$ is given by:
$$
K_{d\omega} = (-a_{46}^2-a_{56}^2)e_1\otimes e^1 +(-a_{46}^2+a_{56}^2) e_2\otimes e^2 +(a_{46}^2+a_{56}^2) e_3\otimes e^3 +(a_{46}^2-a_{56}^2) e_4\otimes e^4 +
$$
$$
+(a_{46}^2-a_{56}^2)e_5\otimes e^5 +(-a_{46}^2+a_{56}^2) e_6\otimes e^6 +2a_{46}a_{56} e_1\otimes e^3 -2a_{46}a_{56} e_3\otimes e^1.
$$
Define the operator $P=K_{d\omega}/|a_{46}^2-a_{56}^2|$.
It defines left-invariant almost para-complex structure, having the property $\omega(PX,PY) = -\omega(X,Y)$. Eigen-subspaces related to eigen-values $(a_{46}^2 -a_{56}^2)/|a_{46}^2 -a_{56}^2|$ and $-(a_{46}^2 -a_{56}^2)/|a_{46}^2 -a_{56}^2|$ of the operator $P$ are generated by the following vectors:
$$
E_1 = \{a_{56} e_1+a_{46}e_3, e_4, e_5 \},\quad E_2 = \{a_{46} e_1+a_{56}e_3, e_2, e_6 \}.
$$
It is easy to see that they are not closed relative to Lie bracket, so $P$ sets non-integrable almost a para-complex structure.

Define the pseudo-Riemannian metric $g(X,Y) = \omega(X,PY)$ of signature (3,3). Direct calculations in the Maple system show that this metric has a nondiagonal Ricci operator and the following scalar curvature
$$
R = \frac{|a_{46}^2-a_{56}^2|}{a_{13}(a_{24}a_{56}-a_{25}a_{46})}.
$$

In particular case, when one of the arguments $a_{46}$ or $a_{56}$ is equal to zero, the situation becomes much simpler. For example, let $a_{56} = 0$.
The property $\omega\wedge d\omega =0$ is fulfilled under the following conditions: $a_{12} = 0$, $a_{15} = 0$,  $a_{35} = 0$. Then 2-form $\omega$ is non-degenerated under the condition $a_{13}a_{25}a_{46}\neq 0$, and we obtain:
$$
\omega =e^1\wedge (a_{13}e^3+a_{14}e^4)+e^2\wedge (a_{23}e^3+a_{24}e^4+a_{25}e^5) +a_{46} e^4\wedge e^6,
$$
$$
d\omega = -a_{46}e^{126} -a_{46} e^{345},
$$
$$
K_{d\omega} = {\rm diag}\{-a_{46}^2, -a_{46}^2, a_{46}^2, a_{46}^2, a_{46}^2, -a_{46}^2\}.
$$
Define the operator $P=K_{d\omega}/a_{46}^2$
It specifies almost para-complex structure, $P^2 = Id$, possessing the property $\omega(PX,PY) = -\omega(X,Y)$.
Eigen-subspaces related to the eigen-values $\pm 1$ of operator $P$ are generated by the following vectors:
$$
E^+ = \{e_3, e_4, e_5 \},\quad E^- = \{e_1, e_2, e_6\}.
$$
It is easy to see that they are not closed relative to Lie bracket, so $P$ sets non-integrable almost a para-complex structure.

The pseudo-Riemannian metric $g(X,Y) =\omega(X,PY)$ of signature (3,3) is given by
$$
g = \left[ \begin {array}{cccccc} 0&0&{\it a_{13}}&{\it a_{14}}&0&0\\
\noalign{\medskip}0&0&{\it a_{23}}&{\it a_{24}}&{\it a_{25}}&0\\
\noalign{\medskip}{\it a_{13}}&{\it a_{23}}&0&0&0&0\\
\noalign{\medskip}{\it a_{14}}&{\it a_{24}}&0&0&0&-{\it a_{46}}\\
\noalign{\medskip}0&{\it a_{25}}&0&0&0&0\\
\noalign{\medskip}0&0&0&-{\it a_{46}}&0&0\end {array} \right]
$$

Direct calculations in the Maple system show that for $a_{14} = a_{24} =0$, this metric has a diagonal Ricci operator with two eigenvalues differing by the sign:
$$
RIC = \frac{a_{46}}{2a_{25}a_{23}}diag\{-1,-1,-1,+1,-1,+1\}.
$$

\subsection{Lie group $G_5$} \label{G_5}
Singular group that does not admit neither symplectic, nor complex structures. Commutation relations:
$$
[e_1,e_2] = e_5, [e_1,e_5] = e_6, [e_3,e_4] = e_6.
$$
Lie algebra g has ideals: $C^1\g = D^1\g = \R\{e_5,e_6\}$, $C^2\g = \R\{e_6\} = \mathcal{Z}$ is the center of Lie algebra. Does admit half-flat structure \cite{Conti}.

Let $\omega = a_{ij}e_i\wedge e_j$ be any 2-form. The Hitchin's operator $K_{d\omega}$ for generic 2-form $\omega$  is given by a quite complicated form. Moreover, $K_{d\omega}^2 = a_{56}^4 Id$.
The $\omega\wedge d\omega = 0$ property is fulfilled under the following conditions:
$$
a_{34} a_{56} + a_{35} a_{46} -a_{36} a_{45} = 0,
$$
$$
a_{12} a_{56} -a_{15}a_{26} + a_{16} a_{25} -a_{23}a_{46} +a_{24} a_{36} -a_{26}a_{34} = 0.
$$
It is easy to see that the 2-form $\omega$ is closed only if
$$
\omega = e^1\wedge (a_{12}e^2+a_{13}e^3+a_{14}e^4+a_{15}e^5) +e^2\wedge (a_{23}e^3+a_{24}e^4 +a_{25}e^5) +a_{34} e^3\wedge e^4.
$$
Such 2-form $\omega$ is non-degenerated.
In order to preserve the non-degeneracy of the $\omega$ and $d\omega$ at the minimal weakening of closedness property of $\omega$, consider the case when $a_{56} \neq 0$.
Then the property $\omega\wedge d\omega =0$ is fulfilled under the following conditions: $a_{34} =0$ and $a_{12} =0$.
The 2-form $\omega$ is non-degenerated under the condition $a_{56}(a_{13}a_{24} -a_{14}a_{23}) \neq 0$, and the following formulas occur:
$$
\omega =e^1\wedge (a_{13}e^3+a_{14}e^4+a_{15}e^5)+e^2\wedge (a_{23}e^3+a_{24}e^4+a_{25}e^5) + a_{56} e^5\wedge e^6,
$$
$$
d\omega = -a_{56} e^{126} +a_{56} e^{345},
$$
$$ 	
K_{d\omega} = a_{56}^2\cdot {\rm diag}\{+1, +1, -1, -1, -1, +1\}.
$$

Define the operator $P= K_{d\omega}/a_{56}^2$.
It defines left-invariant almost para-complex structure $P$, having the property $\omega(PX,PY) = -\omega(X,Y)$.
Eigen-subspaces related to the eigen-values $\pm 1$ of operator $P$ are generated by the following vectors:
$$
E^+ = \{e_1, e_2, e_6 \},\quad E^- = \{e_3, e_4, e_5\}.
$$
It is easy to see that they are not closed relative to Lie bracket, so $P$ sets non-integrable almost a para-complex structure.

The pseudo-Riemannian metric $g(X,Y) = \omega(X,PY)$ of signature (3,3) is given by
$$
g = \left[ \begin {array}{cccccc} 0&0&-{\it a_{13}}&-{\it a_{14}}&-{\it a_{15}}&0
\\ \noalign{\medskip}0&0&-{\it a_{23}}&-{\it a_{24}}&-{\it a_{25}}&0
\\ \noalign{\medskip}-{\it a_{13}}&-{\it a_{23}}&0&0&0&0
\\ \noalign{\medskip}-{\it a_{14}}&-{\it a_{24}}&0&0&0&0
\\ \noalign{\medskip}-{\it a_{15}}&-{\it a_{25}}&0&0&0&{\it a_{56}}
\\ \noalign{\medskip}0&0&0&0&{\it a_{56}}&0\end {array} \right]
$$

Direct calculations in the Maple system show that for $a_{15} = a_{25} =0$, this metric has a diagonal Ricci operator with two eigenvalues differing by the sign:
$$
RIC = \frac{a_{56}}{2a_{13}a_{24}-2a_{14}a_{23}}diag\{-1,-1,-1,-1,+1,+1\}.
$$

\textbf{Conclusions.} Any left-invariant closed 2-form $\omega$ on Lie groups $G_2 - G_5$ is degenerated.
When the closedness requirement of $\omega$ is weakened in order to preserve the non-degeneracy of $\omega$ and $d\omega$ and of the property $\omega\wedge d\omega =0$, the Hitchin's operator $K_{d\omega}$ corresponding to $d\omega$, defines almost para-complex structure $P$.
Pseudo-Riemannian metric $g(X,Y) = \omega(X,PY)$ depending on 5 to 7 arguments is of signature (3,3) and has a diagonal Ricci operator with two eigenvalues differing by the sign.
We have obtained a compatible pair $(\omega,\Omega)$, where $\Omega = d\omega$. Therefore, the properties $d\Omega = 0$ and $\omega\wedge d\omega =0$ were fulfilled in an obvious way.
The para-complex (3,0)-form is given by $\Psi = d\omega +i_\varepsilon \widehat{d\omega}$, where $i_\varepsilon$ is a para-complex unit.
Thus, multiparametric families of almost para-complex half-flat structures were naturally defined on the Lie groups $G_2 - G_5$ and corresponding nilmanifolds. Their structural group is reduced to $SL(3,\R)\subset SO(3,3)$.

\section{Formulas for evaluations}
We now present the formulas which were used for the evaluations (on Maple) of the Nijenhuis tensor and curvature tensor of the associated metrics. Let $e_1,\ldots,e_{2n}$ be a basis of the Lie algebra $\mathfrak g$ and $C_{ij}^k$ a structure constant of the Lie algebra in this base:
\begin{equation}
[e_i,e_j]=\sum_{k=1}^{2n}C_{ij}^{k}e_k,  \label{strukt}
\end{equation}

\textbf{1. Connection components.} These are the components $\Gamma_{ij}^{k}$ in the formula
$\nabla _ {e _ {i}} e_j =\Gamma_{ij}^{k}e_k.$
For left-invariant vector fields we have: $2g ({\nabla}_{X} Y, Z) =g ([X, Y], Z) +g ([Z, X], Y)-g ([Y, Z], X)$. For the basis vectors we have:
$$
2g (\nabla _ {e _ {i}} e_j, e_k) =g ([e_i, e_j], e_k) +g ([e_k, e_i], e_j) +g (e_i, [e_k, e_j]),
$$
$$
2g _ {lk} \Gamma _ {ij} ^ {l} =g _ {pk} C _ {ij} ^ {p} +g _ {pj} C _ {ki} ^ {p} +g _ {ip} C _ {kj} ^ {p},
$$
\begin{equation}
\Gamma_{ij}^{n}=\frac{1}{2}g^{kn}\left(g_{pk}C_{ij}^{p}+g_{pj}C_{ki}^{p} +g_{ip}C_{kj}^{p}\right).
\end{equation}
Maple code:
\begin{maplegroup}
\begin{mapleinput}
\mapleinline{active}{1d}{
Ginv:=inverse(G):
Gamma[i,j,n]:=(1/2)*(sum(Ginv[k,n]*(sum(C[i,j,p]*G[p,k]
+C[k,i,p]*G[p,j]+C[k,j,p]*G[i,p],'p'=1..6)),'k'=1..6)):
}{%
}
\end{mapleinput}
\end{maplegroup}

\vspace{2mm}
\textbf{2. Curvature tensor.}
The formula is: $R(X, Y)Z =\nabla_{X} \nabla_{Y} Z -\nabla_{Y}\nabla_{X}Z -\nabla_{[X,Y]}Z$.
For the basis vectors we have: $R (e_i, e_j) e_k=R _ {ijk} ^s e_s $,
$$
R(e_i,e_j)e_k=\nabla_{e_{i}}\nabla_{e_{j}}e_{k}-\nabla_{e_{j}}\nabla_{e_{i}}e_{k} -\nabla_{[e_{i},e_{j}]}e_{k}.
$$
Therefore:
\begin{equation}
R_{ijk}^{s}=\Gamma_{ip}^{s}\Gamma_{jk}^{p}-\Gamma_{jp}^{s}\Gamma_{ik}^{p} -C_{ij}^{p}\Gamma_{pk}^{s}.
\end{equation}
Maple code:
\begin{maplegroup}
\begin{mapleinput}
\mapleinline{active}{1d}{
Riem[i,j,k,s]:=simplify(sum(Gamma[i,p,s]*Gamma[j,k,p]
-Gamma[j,p,s]*Gamma[i,k,p] -C[i,j,p]*Gamma[p,k,s],'p'=1..6));
}{%
}
\end{mapleinput}
\end{maplegroup}

\vspace{2mm}
\textbf{3. Ricci tensor and scalar curvature.}
The Ricci tensor of a metric $g$ is the tensor $Ric$ obtained by convolution the tensor $R_{ijk}^{s}$ by the first and the upper indices:
\begin{equation}
Ric_{ij}=R_{kij}^{k}.
\end{equation}
Maple code:
\begin{maplegroup}
\begin{mapleinput}
\mapleinline{active}{1d}{
Ric[n,m]:=sum(Riem[i,n,m,i],'i'=1..6);
R:=simplify(sum(sum(Ginv[i,j]*Ric[i,j],'i'=1..6),'j'=1..6));
}{%
}
\end{mapleinput}
\end{maplegroup}

\end{document}